\documentclass[12pt,a4paper,oneside]{article}
\usepackage{amssymb,amsthm,amsmath,fullpage,pictex,dcpic}
\usepackage[english]{babel}

\allowdisplaybreaks[4]
\newenvironment{prova}{\begin{proof}\parindent=0in}{\end{proof}}
\numberwithin{equation}{section}

\newtheorem{prop}{Proposition}[section]
\newtheorem{lemma}[prop]{Lemma}
\newtheorem{df}[prop]{Definition}

\newcommand{\A}{\mathcal{A}}
\newcommand{\B}{\mathcal{B}}
\newcommand{\U}{\mathcal{U}}
\newcommand{\M}{\mathcal{M}}
\newcommand{\HH}{\mathcal{H}}
\newcommand{\C}{\mathbb{C}}
\newcommand{\N}{\mathbb{N}}
\newcommand{\Z}{\mathbb{Z}}
\newcommand{\bn}[1]{\textrm{\footnotesize$\left(\!\!\begin{array}{c}#1 \end{array}\!\!\right)$}}
\newcommand{\inner}[1]{\left<#1\right>}
\newcommand{\ma}[1]{\left(\!\begin{array}{cc}#1 \end{array}\!\right)}

\newcommand{\az}{\triangleright}
\newcommand{\aaz}{\,\textrm{\footnotesize$\blacktriangleright$}_\pm\,}
\newcommand{\za}{\triangleleft}
\newcommand{\bm}{\hspace{1pt}.\hspace{1pt}}
\newcommand{\de}{\mathrm{d}}

\begin{document}

\title{Quantum Groups and Twisted Spectral Triples}

\date{14th February 2007}

\author{~\\[5pt]
        \Large{Francesco D'Andrea} \\ [5pt]
        \normalsize{Scuola Internazionale Superiore di Studi Avanzati,}\\
        \normalsize{Via Beirut 2-4, I-34014, Trieste, Italy}
        \\~\\
}

\maketitle

\bigskip

\thispagestyle{empty}
                                                                               
\begin{abstract}
\noindent
Through the example of the quantum symplectic $4$-sphere, we discuss
how the notion of twisted spectral triple fits into the framework of
quantum homogeneous spaces.
\end{abstract}

\vfill


\pagebreak

\section{Introduction}
In a recent paper~\cite{CM} it is explained how a simple twist in the original
definition of spectral triple~\cite{Connes,CReal} makes it possible to
deal with algebras with no (or few) traces. It was also suggested that
this notion has potential applications to quantum groups and quantum
homogeneous spaces, a domain where to construct a spectral triple
may sometimes be problematic.

We investigate the connection between twisted spectral triples and
quantum homogeneous spaces using as guiding example the algebra
$\A(S^4_q)$ of the quantum symplectic $4$-sphere constructed
in~\cite{LPR}, on which up to now no spectral triples are known.
A possible application of a `twisted' Dirac operator is in the
construction of a differential calculus, which in the particular
case of $S^4_q$ is fundamental for the study of noncommutative
instantons. Another point is the study of the spectral
action~\cite{CC} on $S^4_q$: a problem shared by most quantum
homogeneous spaces is that the axioms for a `real structure'
are fulfilled only modulo an ideal of `infinitesimals', and
we are not able to give meaning to the `adjoint representation'
of $1$-forms (see e.g.~\cite{DLPS,DLSvSV,DDLW,DDL},
and~\cite{DS} for a example which doesn't suffer from this problem).
We'll explain how the notion of real structure can be (trivially)
extended to the case of twisted spectral triples, and then construct
a real structure in the example of $S^4_q$.

We start with a reformulation of the notion of twisted spectral
triple which seems to be appropriate when studying quantum
homogeneous spaces.

\bigskip

We call the data $(\A,\HH,D,K)$ a \emph{twisted spectral triple}
if (i) $\A$ is a complex associative involutive algebra with unity
(for short \emph{$*$-algebra}) represented by bounded operators on a
separable Hilbert space $\HH$, (ii) $D$ is a (unbounded) selfadjoint
operator on $\HH$ with dense domain and compact resolvent, (iii) $K$
is an invertible linear operator on $\mathrm{dom}\hspace{1pt}D$,
(iv) the `$1$-form'
\begin{equation}\label{eq:dea}
\de a:=K^{-1}\bigl(Da-(K^{-1}aK)D\bigr)
\end{equation}
extends to a bounded operator on $\HH$ for any $a\in\A$.
If $K$ is the identity, one gets the original definition of spectral
triple~\cite{Connes}.
If $K$ is bounded the data $(\A,\HH,D',\sigma)$, with $D':=K^{-1}D$ and
$\sigma(a):=K^{-2}aK^2$, is a $\sigma$-spectral triple in the sense
of~\cite{CM}, i.e.~$D'$ has compact resolvent and $D'a-\sigma(a)D'$
is bounded (being equal to $\de a$); notice that not all automorphisms
$\sigma$ are implementable, so the notion of $\sigma$-spectral triple
of~\cite{CM} is more general.
If $K$ and $K^{-1}$ are both bounded, we can `untwist' the Dirac
operator by defining $D''=KD$; since $D''$ has compact resolvent
and $[D'',a]=K^2\de a$ is bounded, the data $(\A,\HH,D'')$ is an
ordinary spectral triple.
As usual, we will refer to $D$ as the `Dirac operator', in analogy with the
commutative situation where spectral triples are canonically associated to
spin structures. Moreover, we'll identify $\A$ with its representation
and omit the representation symbol.

In this paper we construct a triple over the quantum
symplectic $4$-sphere $S^4_q$ which satisfies all the axioms of a
twisted spectral triple, but for the compact resolvent condition.
To compute the spectrum of the Dirac operator is not an easy task,
mainly because $\A(S^4_q)$ has no known symmetries and we cannot
use the powerful tools of representation theory.
We stress that the compact resolvent condition has no role in
the construction of the differential calculus associated to $D$
(not even in the construction of the spectral action),
which we describe in the following.

Given a twisted spectral triple $(\A,\HH,D,K)$, a differential
calculus $(\Omega^\bullet,\de)$ can be constructed as explained in~\cite{CM}.
We define $\Omega^\bullet$ as the $\N$-graded algebra generated by degree $0$
elements $a\in\A$ and degree $1$ elements $\de b$ given by (\ref{eq:dea}),
$b\in\A$, with $\A$-bimodule structure
\begin{equation}\label{eq:Obim}
a\bm\omega=\sigma(a)\omega \;,\qquad
\omega\bm a=\omega a
\end{equation}
where $\sigma(a)=K^{-2}aK^2$ and for all $a\in\A$, $\omega\in\Omega^\bullet$
(on the right hand sides the multiplication in $\B(\HH)$ is understood).
Thanks to the modified bimodule structure, the Leibniz rule is satisfied
$$
\de(ab)=\de a\bm b+a\bm\de b\;,\qquad\forall\;a,b\in\A\;.
$$

The notion of even and real spectral triple can be extended
in a straightforward way to the twisted case.
A (twisted) spectral triple is called \emph{even} if there exists a
$\Z_2$-grading $\gamma$ on $\HH$ (i.e.~a bounded selfadjoint operator
satisfying 
$\gamma^2=1$) such that the Dirac operator is odd and the algebra $\A$ is even:
$$
\gamma D+D\gamma=0\;,\quad\qquad a\gamma=\gamma a \quad\forall\;a\in\A\;.
$$
A \emph{real structure} on a (twisted) spectral triple is
a bounded antilinear operator $J$ on $\HH$ satisfying
\begin{equation}\label{eq:J}
J^2=\pm 1\;,\qquad
JD=\pm DJ\;,
\end{equation}
and such that for all $a,b\in\A$
\begin{equation}\label{eq:real}
[\,a,JbJ^{-1}]=0\;,\qquad
[\,\de a,JbJ^{-1}]=0\;.
\end{equation}
We refer to last equation as the `first order condition'.
If the spectral triple is even, we impose the further condition
$J\gamma=\pm\gamma J$.
The signs `$\pm$' in previous equations are determined by the dimension of
the geometry~\cite{CReal}; a real spectral triple of dimension $4$, for example,
corresponds to the choices $J^2=-1$, $JD=DJ$ and $J\gamma=\gamma J$.

\bigskip

For the reader's ease, we recall also the notions of module algebra
and of crossed product algebra (see e.g.~\cite{KS}), which will be used
throughout the paper.
Let $\A$ be a $*$-algebra and $(\U,\Delta,\epsilon,S)$
a $*$-Hopf algebra. We say that $\A$ is a (left) $\U$-module $*$-algebra
if there is a (left) action `$\az$' of $\U$ on $\A$ satisfying
\begin{equation}\label{eq:covaz}
h\az ab=(h_{(1)}\az a)(h_{(2)}\az b)\;,\qquad
h\az 1=\varepsilon(h)1\;,\qquad
h\az a^*=\{S(h)^*\az a\}^*\;,
\end{equation}
for all $h\in\U$ and $a,b\in\A$. If $\A$ is a (left) $\U$-module
$*$-algebra, the (left) crossed product $\A\rtimes\U$ is defined
as the $*$-algebra generated by $\A$ and $\U$ with crossed
commutation relations
$$
ha=(h_{(1)}\az a)h_{(2)}\;,\quad\forall\;h\in\U,\;a\in\A\,.
$$
Right module algebras and right crossed product algebras
are defined similarly.
As usual we use Sweedler notation for the coproduct,
$\Delta(h)=h_{(1)}\otimes h_{(2)}$.

\bigskip

The plan of the paper is the following. In Section~\ref{sec:due}, we
present the algebra of the quantum symplectic $7$-sphere $S^7_q$
and its symmetry Hopf algebra $U_q(so(5))$. In Section~\ref{sec:tre},
we construct the algebra of the symplectic $4$-sphere $S^4_q$
of~\cite{LPR} as the subalgebra of $\A(S^7_q)$ which is invariant
for the left action of a sub $*$-Hopf algebra of $U_q(so(5))$
isomorphic to $U_q(su(2))$. In Section~\ref{sec:quattro}, we
construct a bounded $*$-representation of $\A(S^4_q)$ on a $\Z_2$-graded
Hilbert space $\HH=\HH_+\oplus\HH_-$ which is a deformation of the space
of $L^2$-spinors on the round $4$-sphere. In Section~\ref{sec:cinque},
we construct a Dirac operator $D$ and prove that, leaving aside
the compact resolvent condition, the data $(\A(S^4_q),\HH,D,K)$
satisfies all other axioms for a twisted spectral triple.
In Section~\ref{sec:sei}, we complete the picture by constructing
a real structure $J$.
The study of the differential calculus on $\A(S^4_q)$ associated
to $D$ is postponed to forthcoming papers. If finite dimensional,
this differential calculus would provide a framework for the
study of $q$-deformations of instantons.

\section{The symplectic $7$-sphere and its symmetries}\label{sec:due}
We now introduce the main characters of this paper, the quantum
universal enveloping algebra $U_q(so(5))$ and the algebra $\A(S^7_q)$
of the symplectic $7$-sphere~\cite{LPR}.

For $0<q<1$, we call $U_q(so(5))$ the `compact' real form of
the Hopf algebra denoted $\breve{U}_q(so(5))$ in~\cite{KS}.
As a $*$-algebra, it is generated by
$\{K_i=K_i^*,K_i^{-1},E_i,F_i:=E_i^*\}_{i=1,2}$ with relations
\begin{equation*}\begin{array}{c}
[K_1,K_2]=0\;\;,\qquad
K_iK_i^{-1}=K_i^{-1}K_i=1\;\;,
\\ \rule{0pt}{20pt}
[E_i,F_j]=\delta_{ij}\frac{K_j^2-K_j^{-2}}{q^j-q^{-j}}\;\;,
\\ \rule{0pt}{20pt}
K_iE_iK_i^{-1}=q^iE_i\;\;,\qquad
K_iE_jK_i^{-1}=q^{-1}E_j\;\;\mathrm{if}\;i\neq j\;,
\end{array}\end{equation*}
plus Serre relations, explicitly, given by
\begin{subequations}\label{eq:Serre}
\begin{align}
E_1E_2^2-(q^2+q^{-2})E_2E_1E_2+E_2^2E_1 &=0\;, \\*
E_1^3E_2-(q^2+1+q^{-2})(E_1^2E_2E_1-E_1E_2E_1^2)-E_2E_1^3 &=0\;.
\end{align}
\end{subequations}
Serre relations can be written in a more compact form by defining
$[a,b]_q:=q^2ab-ba$. Then, (\ref{eq:Serre}) are equivalent to
$$
[E_2,[E_1,E_2]_q]_q=0\;,\qquad
[E_1,[E_1,[E_2,E_1]_q]_q]=0\;.
$$
The Hopf algebra structure $(\Delta,\epsilon,S)$ of $U_q(so(5))$ is given by
\begin{equation*}
\begin{array}{c}
\Delta K_i=K_i\otimes K_i\;\;,\quad \Delta E_i=E_i\otimes K_i+K_i^{-1}\otimes E_i\;\;, \\
\rule{0pt}{3.5ex}\epsilon(K_i)=1\;\;,\quad \epsilon(E_i)=0\;\;,\quad
S(K_i)=K_i^{-1}\;\;,\quad S(E_i)=-q^iE_i\;\;.
\end{array}
\end{equation*}
Notice that the elements $(K_1,K_1^{-1},E_1,F_1)$ generate a sub-Hopf-$*$-algebras
of $U_q(so(5))$, which for obvious reasons we'll denote $U_q(su(2))$ in
the following.

For each non negative $n_1,n_2$ such that $n_2\in\frac{1}{2}\Z$ and $n_2-n_1\in\N$, there
is an irreducible representation of $U_q(so(5))$, whose representation space we denote
$V_{(n_1,n_2)}$ (here $\N$ denotes non-negative integers, i.e.~it includes $0$)
and whose highest weight vector is an eigenvector of $K_1$ and $K_2$
with eigenvalues $q^{n_1}$ and $q^{n_2-n_1}$ respectively. \pagebreak

\noindent
As an example, let us draw the weight diagrams of $V_{(\frac{1}{2},\frac{1}{2})}$,
$V_{(0,1)}$ and $V_{(1,1)}$.

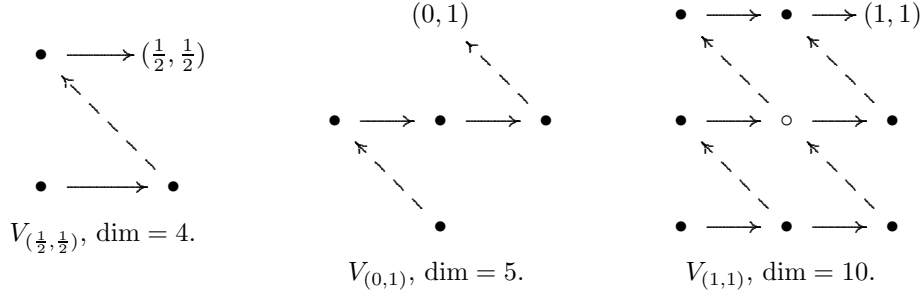
\begin{footnotesize}\begin{center}\begin{tabular}{ccccc}
\begin{tabular}{c}
\begindc{\commdiag}[25]
 \obj(1,1)[A]{$\bullet$}
 \obj(3,1)[B]{$\bullet$}
 \obj(1,3)[C]{$\bullet$}
 \obj(3,3)[D]{$(\frac{1}{2},\frac{1}{2})$}
 \mor{A}{B}{}
 \mor{C}{D}{}
 \mor{B}{C}{}[\atright,\dasharrow]
\enddc
\\ \rule{-10pt}{14pt}$V_{(\frac{1}{2},\frac{1}{2})}$, $\dim=4$.
\end{tabular}
&&~
\begin{tabular}{c}
\begindc{\commdiag}[20]
 \obj(1,3)[A]{$\bullet$}
 \obj(3,3)[B]{$\bullet$}
 \obj(5,3)[C]{$\bullet$}
 \obj(3,1)[D]{$\bullet$}
 \obj(3,5)[E]{$(0,1)$}
 \mor{A}{B}{}
 \mor{B}{C}{}
 \mor{D}{A}{}[\atright,\dasharrow]
 \mor{C}{E}{}[\atright,\dasharrow]
\enddc
\\ \rule{0pt}{14pt}$V_{(0,1)}$, $\dim=5$.
\end{tabular}
&&~
\begin{tabular}{c}
\begindc{\commdiag}[20]
 \obj(3,3)[B]{$\bullet$}
 \obj(5,3)[C]{$\bullet$}
 \obj(7,3)[D]{$\bullet$}
 \obj(3,5)[F]{$\bullet$}
 \obj(5,5)[G]{$\circ$}
 \obj(7,5)[H]{$\bullet$}
 \obj(3,7)[L]{$\bullet$}
 \obj(5,7)[M]{$\bullet$}
 \obj(7,7)[N]{$(1,1)$}
 \mor{B}{C}{}
 \mor{C}{D}{}
 \mor{F}{G}{}
 \mor{G}{H}{}
 \mor{L}{M}{}
 \mor{M}{N}{}
 \mor{C}{F}{}[\atright,\dasharrow]
 \mor{D}{G}{}[\atright,\dasharrow]
 \mor{G}{L}{}[\atright,\dasharrow]
 \mor{H}{M}{}[\atright,\dasharrow]
\enddc
\\ \rule{-10pt}{14pt}$V_{(1,1)}$, $\dim=10$.
\end{tabular}
\end{tabular}\end{center}\end{footnotesize}

\noindent A solid arrow indicates points that can be joined by applying $E_1$ (the reverse
arrow corresponds to $F_1$), a dashed arrow indicates points that can be joined by applying
$E_2$ (the reverse arrow corresponds to $F_2$).
In each diagram, the highest weight vector is denoted by its weight $(n_1,n_2)$. A bullet
indicates a weight with multiplicity $1$, an empty circle a weight with multiplicity $2$.

The Hopf algebra $U_q(so(5))$ describes the symmetries of an
algebra which is a deformation of the algebra of polynomial
functions on the $7$-sphere $S^7$.

\begin{df}\label{def:Sseven}
We call symplectic $7$-sphere the `virtual space' underlying the $*$-algebra
$\A(S^7_q)$ with generators $z_i,z_i^*$ ($i=1,\ldots,4$), commutation relations
\begin{align*}
z_1z_2&=q^{-1}z_2z_1\;, &
z_2z_4&=q^{-1}z_4z_2\;, &
z_2z_4^*&=q^{-1}z_4^*z_2 \;, \\
z_1z_3&=q^{-1}z_3z_1\;, &
z_3z_4&=q^{-1}z_4z_3\;, &
z_3z_4^*&=q^{-1}z_4^*z_3 \;, \\
z_1z_4&=q^{-2}z_4z_1\;, &
z_1z_4^*&=q^{-2}z_4^*z_1\;, &
z_2z_3^*&=q^{-2}z_3^*z_2\;, \\
z_2z_3-q^2z_3z_2&=q(1-q^2)z_1z_4 \;, &
z_1z_2^*-q^{-1}z_2^*z_1&=(q-q^{-1})z_4^*z_3 \;,&
z_1z_3^*-q^{-1}z_3^*z_1&=q(1-q^2)z_4^*z_2 \;, \\
[z_4^*,z_4]&=0 \;, &
[z_2^*,z_2]&=(1-q^2)z_4z_4^* \;, &
[z_3^*,z_3]&=z_2^*z_2-q^4z_2z_2^* \;,
\end{align*}
and
$$
z_1z_1^*+z_2z_2^*+z_3z_3^*+z_4z_4^*=
z_1^*z_1+q^6z_2^*z_2+q^2z_3^*z_3+q^8z_4^*z_4=1 \;.
$$
\end{df}

\smallskip

\noindent
The generators used in~\cite{LPR}, which we denote $x'_i$, are related
to ours by the equations $x'_1=q^4z_4$, $x'_2=q^3z_2$, $x'_3=-qz_3$ and $x'_4=z_1$.

\begin{prop}
The algebra $\A(S^7_q)$ is an $U_q(so(5))$-module $*$-algebra
for the action defined on generators by
\begin{align*}
\rule{0pt}{16pt}
K_1\az z_1&=q^{1/2}z_1\;, &
K_1\az z_2&=q^{1/2}z_2\;, &
K_1\az z_3&=q^{-1/2}z_3\;, &
K_1\az z_4&=q^{-1/2}z_4\;, \\
\rule{0pt}{16pt}
K_2\az z_1&=z_1\;, &
K_2\az z_2&=q^{-1}z_2\;, &
K_2\az z_3&=qz_3\;, &
K_2\az z_4&=z_4\;, \\
\rule{0pt}{16pt}
E_1\az z_1&=0\;, &
E_1\az z_2&=0\;, &
E_1\az z_3&=z_1\;, &
E_1\az z_4&=z_2\;, \\
\rule{0pt}{16pt}
E_2\az z_1&=0\;, &
E_2\az z_2&=z_3\;, &
E_2\az z_3&=0\;, &
E_2\az z_4&=0\;, \\
\rule{0pt}{16pt}
F_1\az z_1&=z_3\;, &
F_1\az z_2&=z_4\;, &
F_1\az z_3&=0\;, &
F_1\az z_4&=0\;, \\
\rule{0pt}{16pt}
F_2\az z_1&=0\;, &
F_2\az z_2&=0\;, &
F_2\az z_3&=z_2\;, &
F_2\az z_4&=0\;.
\end{align*}
Let $U_q(u(1))$ be the Hopf $*$-algebra generated by $K_1,K_1^{-1}$.
Then $\A(S^7_q)$ is also a right $U_q(u(1))$-module $*$-algebra
for the action defined by $z_i\za K_1=q^{1/2}z_i$. Left and
right actions commute.
\end{prop}

\begin{prova}
We consider first the free $*$-algebra generated by $\{z_i,z_i^*\}$.
The elements $z_i$ carry the fundamental representation
$V_{(\frac{1}{2},\frac{1}{2})}$ of $U_q(so(5))$, while the action
is extended to $z_i^*$ by compatibility with the involution using
the rule $h\az a^*=\{S(h)^*\az a\}^*$. Thus the free $*$-algebra
is trivially an $U_q(so(5))$-module $*$-algebra.

Let $V:=V_{(\frac{1}{2},\frac{1}{2})}\oplus V_{(\frac{1}{2},\frac{1}{2})}^*$.
Inside the decomposition of Hopf tensor product $V\otimes V$ we consider
the following vectors (tensor product symbol implied)
\begin{align*}
v_1 &:=z_1z_4^*-q^{-2}z_4^*z_1 \;,&
v_5 &:=q^2z_1z_4-q^{-2}z_4z_1+q^{-1}z_2z_3-qz_3z_2 \;,\\
v_2 &:=z_1z_3-q^{-1}z_3z_1 \;,&
v_6 &:=v_5^* \;,\\
v_3 &:=(z_2z_4-q^{-1}z_4z_2)^* \;,&
v_7 &:=z_1^*z_1+q^6z_2^*z_2+q^2z_3^*z_3+q^8z_4^*z_4-1 \;,\\
v_4 &:=z_1z_2^*+z_3z_4^*-q^{-1}(z_2^*z_1+q^2z_4^*z_3) \;,&
v_8 &:=z_1z_1^*+z_2z_2^*+z_3z_3^*+z_4z_4^*-1 \;.
\end{align*}
As one check by direct computation, $v_i$'s are eigenvectors of
$K_i$'s and annihilated by both $E_1$ and $E_2$, hence they are
highest weight vectors of irreducible representations of $U_q(so(5))$
inside $V\otimes V$.
In particular, $v_1$ has weight $(1,1)$, $\{v_2,v_3,v_4\}$
have weight $(0,1)$, $\{v_5,v_6,v_7,v_8\}$ have weight $(0,0)$.
By applying $F_1$ and $F_2$ to these highest weight vectors
one proves that a linear basis for the \emph{real}
$29$-dimensional representation
$$
4V_{(0,0)}\oplus 3V_{(0,1)}\oplus V_{(1,1)}\subset V\otimes V
$$
is given just by the degree $\leq 2$ polynomials appearing
in Proposition \ref{def:Sseven} and by their conjugated.
Thus the ideal they generate is $U_q(so(5))$-invariant
and $\A(S^7_q)$, being the quotient of an $U_q(so(5))$-module
$*$-algebra by a two-sided invariant $*$-ideal, is itself
an $U_q(so(5))$-module $*$-algebra.

The remaining part of the Proposition is trivial.
\end{prova}

\smallskip

\noindent
From the general theory of compact matrix quantum groups
(cf.~Section 11 of \cite{KS}) we know that there
exists a (unique) positive faithful $U_q(so(5))$-invariant
linear functional $\varphi:\A(S^7_q)\to\C$. This functional
comes from the Haar state of the Hopf algebra $\A(Sp_q(2))$
dual to $U_q(so(5))$, and satisfies
\begin{equation}\label{eq:twist}
\varphi(ab)=\varphi\bigl(b\,\kappa(a)\bigr)\;,
\end{equation}
for all $a,b\in\A(S^7_q)$, where $\kappa:\A(S^7_q)\to\A(S^7_q)$
is called the `modular automorphism' and is given by
(cf.~Section 11.3.4 of \cite{KS})
\begin{equation}\label{eq:mod}
\kappa(a)=K_1^8K_2^6\az a\za K_1^8 \;.
\end{equation}

\section{The quantum symplectic $4$-sphere}\label{sec:tre}
Consider the algebra
\begin{equation}\label{eq:Squattro}
\A(S^4_q):=\big\{ a\in\A(S^7_q) \,\big|\, h\az a=\epsilon(h)a\;\forall\;h\in U_q(su(2)) \big\}
\end{equation}
where $U_q(su(2))$ is the sub-Hopf-$*$-algebra of $U_q(so(5))$ generated by
$(K_1,K_1^{-1},E_1,F_1)$. We now show that this is just the algebra of `functions'
of the $4$-sphere constructed in \cite{LPR}.

Dually to the left action of $U_q(su(2))$ we can define a
right coaction of $\A(SU_q(2))$ such that $\A(S^7_q)$ is a
right comodule $*$-algebra.
Recall that $\A(SU_q(2))$ is the $*$-algebra generated by
$\alpha,\beta$ and the adjoints with relations
\begin{equation*}
\beta\alpha=q\alpha\beta\;\;,\quad
\beta^*\alpha=q\alpha\beta^*\;\;,\quad
[\beta,\beta^*]=0\;\;,\quad
\alpha^*\alpha+q^2\beta^*\beta=1\;\;,\quad
\alpha\alpha^*+\beta\beta^*=1\;\;,
\end{equation*}
with Hopf algebra structure
\begin{gather*}
\Delta\ma{\alpha & \beta \\ -q\beta^* & \alpha^*}
=\ma{\alpha & \beta \\ -q\beta^* & \alpha^*}\,\dot{\otimes}\,
\ma{\alpha & \beta \\ -q\beta^* & \alpha^*}\;\;, \\
\rule{0pt}{25pt}
\epsilon\ma{\alpha & \beta \\ -q\beta^* & \alpha^*}=1\;\;,\qquad
S\ma{\alpha & \beta \\ -q\beta^* & \alpha^*}=
\ma{\alpha^* & -q\beta \\ \beta^* & \alpha}\;\;,
\end{gather*}
and with obvious $*$-structure. We use the same notation of~\cite{DLSvSV},
but for greek letters instead of latin ones. The dotted tensor product
is defined as $(A\,\,\dot{\otimes}\,B)_{ij}=\sum_kA_{ik}\otimes B_{kj}$.

The coaction $\Delta_R(a)=a_{(0)}\otimes a_{(1)}$ dual to the right action
of $U_q(su(2))$ is determined by:
$$
h\az a=:a_{(0)}\inner{h,a_{(1)}}\;,
$$
with $\inner{\,\,,\,}$ the dual pairing between $U_q(su(2))$ and $\A(SU_q(2))$,
given by (cf.~Section 4.4.1 of \cite{KS})
$$
\inner{K_1,\alpha}=q^{1/2}\;,\qquad\inner{K_1,\alpha^*}=q^{-1/2}\;,\qquad
\inner{E_1,\beta}=\inner{F_1,-q\beta^*}=1\;, \\
$$
while $\inner{h,a}=0$ for any other pair of generators.
Using these we compute the coaction dual to the left action of $U_q(su(2))$,
and get
\begin{align*}
\Delta_R(z_1)   &=z_1\otimes\alpha-qz_3\otimes\beta^* \;, \\
\Delta_R(z_2)   &=z_2\otimes\alpha-qz_4\otimes\beta^* \;, \\
\Delta_R(z_3)   &=z_3\otimes\alpha^*+z_1\otimes\beta \;, \\
\Delta_R(z_4)   &=z_4\otimes\alpha^*+z_2\otimes\beta \;.
\end{align*}
Let $\Psi$ be the matrix
\begin{equation}\label{eq:Psi}
\Psi=\left(\begin{array}{rr}
  -qz_3^* &    z_1^* \\
      z_1 &    z_3   \\
     qz_2 &   qz_4   \\
-q^3z_4^* & q^2z_2^*
\end{array}\right) \;.
\end{equation}
The coaction can be encoded in the compact formula
\begin{equation}\label{eq:Coa}
\Delta_R(\Psi)=\Psi\,\dot{\otimes}\,\ma{\alpha & \beta \\ -q\beta^* & \alpha^*}\;.
\end{equation}
Since $\Psi^\dag\Psi=1$, the matrix $P:=\Psi\Psi^\dag$ is a projection and by
construction its elements are $\A(SU_q(2))$-coinvariant.
By duality between the action of $U_q(su(2))$ and the coaction of $\A(SU_q(2))$,
the algebra of coinvariant is just $\A(S^4_q)$, which then coincides
with the algebra called symplectic $4$-sphere in \cite{LPR}. In particular,
in \cite{LPR} it was proved that $\A(S^4_q)$ is generated by the matrix
elements $P_{ij}$ of $P$, with relations and $*$-structure $P=P^*=P^2$.

We choose the generators
$$
x_0:=z_2z_2^*+z_4z_4^*\;,\qquad
x_1:=q(z_1z_2^*+z_3z_4^*)\;,\qquad
x_2:=z_2z_3-qz_1z_4\;,
$$
and notice that the projection $P$ becomes
\begin{equation}\label{eq:P}
P=\left(\begin{array}{cccc}
    1-q^6x_0 & 0 & -q^2x_2^* & q^2x_1^* \\
    0 & 1-x_0 & x_1 & x_2 \\
    -q^2x_2 & x_1^* & q^2x_0 & 0 \\
    q^2x_1 & x_2^* & 0 & q^4x_0
\end{array}\right) \;.
\end{equation}
With this expression, we compute the relations among the $x_i$'s, which
we summarize in the next proposition.

\begin{prop}
The $*$-algebra $\A(S^4_q)$ is generated by $x_0=x_0^*$, $x_i$ and $x_i^*$
($i=1,2$) with relations
$$
x_0x_i=q^{2i}x_ix_0 \;,\quad
x_1x_2=x_2x_1 \;,\quad
x_1^*x_2=q^4x_2x_1^* \;,\quad
x_i^*x_i-q^4x_ix_i^*=(1-q^{2i})(q^2x_0)^i \;,
$$
and
$$
x_0^2+x_1x_1^*+x_2x_2^*=x_0\;.
$$
\end{prop}

\medskip

\noindent The generators used in~\cite{LPR} are given by
$t:=q^4x_0$, $a:=x_1^*$ and $b:=q^2x_2$.

\section{The modules of chiral spinors}\label{sec:quattro}
The left regular representation of a compact matrix quantum group is bounded
with respect to the inner product induced by the Haar state $\varphi$.
Being $\A(S^7_q)\subset\A(Sp_q(2))$, if we define the inner product
of two elements $v=(v_1,v_2)$ and $w=(w_1,w_2)$ of $\A(S^7_q)^2$ as
\begin{equation}\label{eq:in}
\inner{v,w}:=\varphi(v_1^*w_1)+\varphi(v_2^*w_2)\;,
\end{equation}
then the representation of $\A(S^7_q)$ by left multiplication on this inner
product space is a bounded $*$-representation. We now define two
subspaces of $\A(S^7_q)^2$ which are invariant when multiplied by elements
in the subalgebra $\A(S^4_q)\subset\A(S^7_q)$.

Let $\sigma_+:U_q(so(2))\to\mathrm{Mat}_4(\C)$ be following representation
$$
\sigma_+(K_1)=\ma{q^{1/2} & 0 \\ 0 & q^{-1/2}}\;,\qquad
\sigma_+(E_1)=\ma{0 & 1 \\ 0 & 0}\;,\qquad
\sigma_+(F_1)=\ma{0 & 0 \\ 1 & 0}\;,
$$
where as before $U_q(su(2))$ denotes the $*$-Hopf algebra
generated by $(K_1,K_1^{-1},E_1,F_1)$.
A second (unitary equivalent) representation is given by
$$
\sigma_-(h)=\ma{0 & -1 \\ 1 & 0}\sigma_+(h)\ma{0 & 1 \\ -1 & 0}\;.
$$
for all $h\in U_q(su(2))$.

If we write $v=(v_1,v_2)\in\A(S^7_q)^2$ as a rwo vector,
two actions of $h\in U_q(su(2))$ on $\A(S^7_q)^2$ can be
defined through the formula
\begin{equation}\label{eq:blackaz}
h\aaz v:=(h_{(1)}\az v)\sigma_\pm(S(h_{(2)}))\;,
\end{equation}
where on the right row by column multiplication is understood.
With these actions, we define two subspaces $\M_\pm$ of $\A(S^7_q)^2$
as follows
\begin{equation}\label{eq:emme}
\M_\pm:=\big\{ v\in\A(S^7_q)^2 \,\big|\,
h\aaz v=\epsilon(h)v\;\forall\;h\in U_q(su(2))\big\}\;,
\end{equation}
which in particular means that
\begin{equation}\label{eq:condP}
K_1\az (v_1,v_2)=(q^{\frac{1}{2}}v_1,q^{-\frac{1}{2}}v_2)\;,\qquad
E_1\az (v_1,v_2)=(0,qv_1)\;,\qquad
F_1\az (v_1,v_2)=(q^{-1}v_2,0)\;,
\end{equation}
for all $(v_1,v_2)\in\M_+$ and
\begin{equation}\label{eq:condM}
K_1\az (w_1,w_2)=(q^{-\frac{1}{2}}w_1,q^{\frac{1}{2}}w_2)\,,\;\;\;
E_1\az (w_1,w_2)=(-qw_2,0)\,,\;\;\;
F_1\az (w_1,w_2)=(0,-q^{-1}w_1)\,,
\end{equation}
for all $(w_1,w_2)\in\M_-$.
Conditions (\ref{eq:condP}) are necessary and sufficient for a vector
$(v_1,v_2)$ to be an element of $\M_+$, and conditions (\ref{eq:condM})
are necessary and sufficient for a vector $(w_1,w_2)$ to be an element
of $\M_-$.

\begin{lemma}
The linear spaces $\M_+$ and $\M_-$ are orthogonal.
\end{lemma}
\begin{prova}
From (\ref{eq:condP}) and (\ref{eq:condM}) we get
$$
K_1\az v_1^*w_1=q^{-1}v_1^*w_1\;,\qquad
K_1\az v_2^*w_2=qv_2^*w_2\;.
$$
for all $v\in\M_+$ and $w\in\M_-$.
Applying the Haar functional to both sides of the equations and
using its invariance we get
$$
\varphi(v_1^*w_1)=q^{-1}\varphi(v_1^*w_1)\;,\qquad
\varphi(v_2^*w_2)=q\varphi(v_2^*w_2)\;,
$$
which imply $\varphi(v_1^*w_1)=\varphi(v_2^*w_2)=0$ and
then $\inner{v,w}=0$ for all $v\in\M_+$ and $w\in\M_-$.
\end{prova}

\begin{lemma}
The linear spaces $\M_\pm$ are $\A(S^4_q)$-bimodules.
\end{lemma}
\begin{prova}
By (\ref{eq:Squattro}) we have
\begin{align*}
h\aaz (av) &=(h_{(1)}\az a)(h_{(2)}\aaz v)=a\{\epsilon(h_{(1)})h_{(2)}\aaz v\}=a(h\aaz v) \\
h\aaz (va) &=h_{(1)}\az (va)\sigma_\pm(S(h_{(2)}))=
   \{h_{(1)}\az v\sigma_\pm(S(h_{(3)}))\}(h_{(2)}\az a)=(h\aaz v)a
\end{align*}
for all $h\in U_q(su(2))$, $a\in\A(S^4_q)$ and $v\in\A(S^7_q)^2$.
In particular, if $v$ is an invariant element for the action $\aaz$,
$av$ and $va$ are invariant elements too. Hence, $\M_\pm$ are invariant
subspaces of $\A(S^7_q)^2$ with respect to left/right multiplication by
$\A(S^4_q)$.
\end{prova}

We thus have two bounded $*$-representations of $\A(S^4_q)$, denoted
$\tilde{\pi}_\pm$ and given by
\begin{subequations}\label{eq:repsPM}
\begin{align}
\tilde{\pi}_+(a)v &=a\cdot v\;,\qquad a\in\A(S^4_q),\;v\in\M_+\,, \\
\tilde{\pi}_-(a)w &=a\cdot w\;,\qquad a\in\A(S^4_q),\;w\in\M_-\,.
\end{align}
\end{subequations}
We denote $(\HH_\pm,\pi_\pm)$ the Hilbert space completion of
these representations, and call \emph{spinorial representation}
$(\HH,\pi)$ their direct sum, i.e.~$\pi=\pi_+\oplus\pi_-$ and
$\HH=\HH_+\oplus\HH_-$.
This representation is even with respect to the natural grading $\gamma$
on $\HH$ which assigns degree $+1$ (resp.~$-1$) to the Hilbert space
$\HH_+$ (resp.~$\HH_-$). The data $(\A(S^4_q),\HH,\gamma)$ is the building
block of a twisted spectral triple, to which in the next sections
we will add a Dirac operator $D$ and a real structure $J$.

Note that since $\M_+$ and $\M_-$ are orthogonal,
their direct sum is just their linear span.

\bigskip

Usually, an additional requirement for a spectral triple
is that the domain of the Dirac operator contains a dense subspace
of $\HH$ which is a finitely generated projective (left) $\A$-module
(this is called `finiteness axiom' in Section 10.5 of~\cite{GVF}).
With this in mind, we now prove that $\M_+$ is finitely generated
and projective both as left and right $\A(S^4_q)$-module.

\begin{prop}\label{prop:left}
$\M_+$ is isomorphic to $\A(S^4_q)^4P$ as a left $\A(S^4_q)$-module,
with $P=\Psi\Psi^\dag$ and $\Psi$ the matrix defined in Equation
(\ref{eq:Psi}).
\end{prop}
\begin{prova}
A linear map $\rho:\M_+\to\A(S^4_q)^4P$ is defined by
$$
\rho(v_1,v_2)=(v_1,v_2)\cdot
\textrm{\footnotesize$\left(\!\!\begin{array}{cc}
q & 0 \\ 0 & 1\end{array}\!\!\right)$}\Psi^\dag \;,
$$
where matrix multiplication on the right hand side is
understood, and writing elements of $\A(S^4_q)^4P$ as
row vectors.
Since $\Psi^\dag P=\Psi^\dag$ it is clear that the image of
$\rho$ is in $\A(S^4_q)^4P$. By construction $\rho$ is a
left $\A(S^4_q)$-module map.

A second left $\A(S^4_q)$-module map $\rho^{-1}:\A(S^4_q)^4P\to\M_+$ is
given by
$$
\rho^{-1}(a_1,a_2,a_3,a_4):=(a_1,a_2,a_3,a_4)\cdot\Psi
\textrm{\footnotesize$\left(\!\!\begin{array}{cc}
q^{-1} & 0 \\ 0 & 1\end{array}\!\!\right)$}\;.
$$
From (\ref{eq:Psi}), the invariance of $a_i$'s
and the explicit expression of the action of $(K_1,E_1,F_1)$
on $\A(S^7_q)$, one proves that $\rho^{-1}(a_1,a_2,a_3,a_4)$
satisfies (\ref{eq:condP}). Hence the image of $\rho^{-1}$ is
in $\M_+$.

Since $\Psi^\dag\Psi=1$ and right multiplication for
$\Psi\Psi^\dag$ is the identity operator on $\A(S^4_q)P$,
the maps $\rho$ and $\rho^{-1}$ are one the inverse of the
other (as the notation suggests) and so $\rho$ is a bijective
left $\A(S^4_q)$-module map,
i.e.~an isomorphism of left $\A(S^4_q)$-modules.
\end{prova}

\begin{prop}\label{prop:right}
$\M_+$ is isomorphic to $P\A(S^4_q)^4$ as a right $\A(S^4_q)$-module.
\end{prop}
\begin{prova}
We write elements of the right projective module $P\A(S^4_q)^4$
as column vectors. The proof is similar to the proof of
Proposition \ref{prop:left}, with the only difference that
now the maps realizing the isomorphism $\M_+\simeq P\A(S^4_q)^4$,
which we denote again $\rho:\M_+\to P\A(S^4_q)^4$ and
$\rho^{-1}:P\A(S^4_q)^4\to\M_+$, are given by
$$
\rho(v_1,v_2)=\Psi\textrm{\footnotesize$\left(\!\!\begin{array}{cc}
0 & 1 \\ -1 & 0\end{array}\!\!\right)$}(v_1,v_2)^t\;,\qquad
\rho^{-1}\bn{a_1 \\ a_2 \\ a_3 \\ a_4}=\left\{
\textrm{\footnotesize$\left(\!\!\begin{array}{cc}
0 & -1 \\ 1 & 0\end{array}\!\!\!\right)$}\Psi^\dag
\bn{a_1 \\ a_2 \\ a_3 \\ a_4}\right\}^t\;.
$$
Clearly $\rho$ and $\rho^{-1}$ are right $\A(S^4_q)$-linear,
are one the inverse of the other, the image of $\rho$ is
in $P\A(S^4_q)^4$ simply because $P\Psi=\Psi$, and
the image of $\rho^{-1}$ is $\M_+$ since $\rho^{-1}(a)$
satisfies (\ref{eq:condP}) for all $a_i\in\A(S^4_q)$.
\end{prova}

\section{The Dirac operator}\label{sec:cinque}
Let $v=(v_1,v_2)\in\A(S^7_q)^2$ and consider (for $i=1,2,3$) the linear maps
$$
v\mapsto X^iv:=\sum\nolimits_{j=1,2}
\bigl(X^i_{1,j}\az v_j,X^i_{2,j}\az v_j\bigr)
$$
with $((X^i_{jk}))$ the matrices
\begin{align*}
((X^1_{jk})) &=\ma{
      q[2]E_2 & q[E_2,E_1]_q \\
 q^{-1}[E_2,E_1]_q & -[E_1,[E_2,E_1]_q] } \;, \\
((X^2_{jk})) &=\ma{
-[F_1,[F_1,F_2]_q] & q[F_1,F_2]_q \\
 q^{-1}[F_1,F_2]_q & -q[2]F_2  } \;, \\
((X^3_{jk})) &=K_2^{-1}\ma{
   0 & -1 \\
   q & 0 } \;.
\end{align*}

\begin{lemma}\label{lemma:Xmap}
The operator $X^i$ maps $\M_+$ into $\M_-$, for all $i=1,2,3$.
\end{lemma}
\begin{prova}
Let $v$ be a vector in $\M_+$, thus satisfying (\ref{eq:condP}),
We want to prove that the vector $(w_1,w_2):=X^i(v_1,v_2)$ satisfies the
conditions (\ref{eq:condM}) defining $\M_-$. If $i=3$ the check is
trivial. Let us focus on the cases $i=1,2$.

From the defining relations of $U_q(so(5))$ we derive the following
commutation rules
\begin{align*}
K_1A &=q^{-1}AK_1 &&\mathrm{for}\;A=E_2,[E_1,[E_2,E_1]_q] \;,\\
K_1A &=qAK_1      &&\mathrm{for}\;A=F_2,[F_1,[F_1,F_2]_q] \;,\\
K_1A &=AK_1       &&\mathrm{for}\;A=[E_2,E_1]_q,[F_1,F_2]_q \;,\\
[F_1,[E_2,E_1]_q] &=-[2]E_2K_1^2  \;,\\
[E_1,[F_1,F_2]_q] &=q^2[2]F_2K_1^2 \;,\\
E_1[E_1,[E_2,E_1]_q] &=q^{-2}[E_1,[E_2,E_1]_q]E_1 
 && (\textrm{Serre relation}) \;,\\
F_1[F_1,[F_1,F_2]_q] &=q^2[F_1,[F_1,F_2]_q]F_1 
 && (\textrm{Serre relation}) \;.
\end{align*}
The first three equations imply $K_1\az (w_1,w_2)=(q^{-1/2}w_1,q^{1/2}w_2)$,
which is the first condition of (\ref{eq:condM}), while the remaining equations
gives us for $i=1$
\begin{align*}
F_1\az w_1 &=q[2]E_2F_1\az v_1+q[F_1,[E_2,E_1]_q]\az v_2
 && (\mathrm{since}\;F_1\az v_2=0) \\
           &=q[2]E_2F_1\az v_1-q[2]E_2K_1^2\az v_2 \\
           &=q[2]E_2\az (F_1\az v_1-q^{-1}v_2)=0 \;, \\
E_1\az w_1 &=q[2]E_1E_2\az v_1+qE_1[E_2,E_1]_q\az v_2 \\
           &=-[E_2,E_1]_q\az v_1+q[E_1,[E_2,E_1]_q]\az v_2
 && (\mathrm{since}\;E_1\az v_j=q\delta_{j2}v_2) \\
           &=-qw_2 \;, \\
E_1\az w_2 &=q^{-1}E_1[E_2,E_1]_q\az v_1-E_1[E_1,[E_2,E_1]_q]\az v_2 \\
           &=q^{-1}E_1[E_2,E_1]_q\az v_1-q^{-2}[E_1,[E_2,E_1]_q]E_1\az v_2 \\
           &=q^{-1}[E_2,E_1]_qE_1\az v_1=0 \;,
\end{align*}
and for $i=2$
\begin{align*}
F_1\az w_1 &=-F_1[F_1,[F_1,F_2]_q]\az v_1+qF_1[F_1,F_2]_q\az v_2 \\
           &=-q^2[F_1,[F_1,F_2]_q]F_1\az v_1+qF_1[F_1,F_2]_q\az v_2 \\
           &=q[F_1,F_2]_qF_1\az v_2=0 \;, \\
F_1\az w_2 &=q^{-1}F_1[F_1,F_2]_q\az v_1-q[2]F_1F_2\az v_2 \\
           &=q^{-1}[F_1,[F_1,F_2]_q]\az v_1-[F_1,F_2]_q\az v_2
 && (\mathrm{since}\;F_1\az v_j=q^{-1}\delta_{j1}v_2) \\
           &=-q^{-1}w_1 \;, \\
E_1\az w_2 &=q^{-1}[E_1,[F_1,F_2]_q]\az v_1-q[2]F_2E_1\az v_2
 && (\mathrm{since}\;E_1\az v_1=0) \\
           &=q[2]F_2\az (qv_1-E_1\az v_2)=0 \;.
\end{align*}
These equations implies
\begin{align*}
qF_1\az w_2&=-F_1E_1\az w_1=[E_1,F_1]\az w_1=\tfrac{K_1^2-K_1^{-2}}{q-q^{-1}}\az w_1=-w_1
   && \mathrm{if}\; i=1 \;,\\
q^{-1}E_1\az w_1&=-E_1F_1\az w_2=-[E_1,F_1]\az w_2=-\tfrac{K_1^2-K_1^{-2}}{q-q^{-1}}\az w_2=-w_2
   && \mathrm{if}\; i=2 \;.
\end{align*}
Thus, $(w_1,w_2)$ satisfies all the conditions in Equation (\ref{eq:condM}),
and this concludes the proof.
\end{prova}

\begin{lemma}\label{lemma:XmapS}
The operator $(X^i)^*$ maps $\M_-$ into $\M_+$, for all $i=1,2,3$.
\end{lemma}
\begin{prova}
Similar to the proof of Lemma \ref{lemma:Xmap}.
\end{prova}

\begin{lemma}\label{lemma:Xbound}
The operators
\begin{equation}\label{eq:ltp}
K_2^{-1}\bigl(X^ia-(K_2^{-1}\az a)X^i\bigr)\qquad\textrm{and}\qquad
\bigl(X^ia-(K_2^{-1}\az a)X^i\bigr)K_2^{-1}
\end{equation}
are bounded on $\M_+$, for all $a\in\A(S^4_q)$ and for all $i=1,2,3$.
\end{lemma}
\begin{prova}
The case $i=3$ is trivial, being $X^3a=(K_2^{-1}\az a)X^3$ by
covariance of the action.
For the remaining cases, we proceed as follows.
Firstly, by direct computation one proves
\begin{subequations}\label{eq:subDelta}
\begin{align}
\Delta(E_2)&=K_2^{-1}\otimes E_2+E_2\otimes K_2 \\
\Delta(F_2)&=K_2^{-1}\otimes F_2+F_2\otimes K_2 \\
\Delta([E_2,E_1]_q)&\sim K_2^{-1}\otimes [E_2,E_1]_q
             +[E_2,E_1]_q\otimes K_1K_2 \label{eq:C} \\
\Delta([F_1,F_2]_q)&\sim K_2^{-1}\otimes [F_1,F_2]_q
             +[F_1,F_2]_q\otimes K_1K_2 \\
\Delta([E_1,[E_2,E_1]_q])&\sim K_2^{-1}\otimes [E_1,[E_2,E_1]_q]
             +[E_1,[E_2,E_1]_q]\otimes K_1^2K_2 \\
\Delta([F_1,[F_1,F_2]_q])&\sim K_2^{-1}\otimes [F_1,[F_1,F_2]_q]
             +[F_1,[F_1,F_2]_q]\otimes K_1^2K_2
\end{align}
\end{subequations}
where the notation $A\sim B$ means that the difference
$A-B$ acts trivially on $\A(S^4_q)\otimes\A(S^7_q)$.
For example
$$
\Delta([E_2,E_1]_q)=
(K_1K_2)^{-1}\otimes [E_2,E_1]_q
+[E_2,E_1]_q\otimes K_1K_2
+(q^2-q^{-2})K_2^{-1}E_1\otimes E_2K_1 \;,
$$
hence using $E_1\otimes 1\sim 0$ and $K_1\otimes 1\sim 1$ we get
the expression in Equation (\ref{eq:C}).

Then, we observe that each matrix element $X^i_{jk}$ of $X^i$
($i,j,k=1,2$) is proportional to one of the elements whose coproduct
has been computed in (\ref{eq:subDelta}). Thus,
$$
\Delta(X^i_{jk})\sim K_2^{-1}\otimes X^i_{jk}+X^i_{jk}\otimes K_1^{n_{ijk}}K_2
$$
for all $i,j,k$ and for a suitable $n_{ijk}\geq 0$.
This equation together with (\ref{eq:covaz}) implies
\begin{subequations}\label{eq:Bcheck}
\begin{align}
K_2^{-1}\bigl(X^i_{jk}a-(K_2^{-1}\az a)X^i_{jk}\bigr)&=(K_2^{-1}X^i_{jk}\az a)K_1^{n_{ijk}}\;,\\
\bigl(X^i_{jk}a-(K_2^{-1}\az a)X^i_{jk}\bigr)K_2^{-1}&=(X^i_{jk}\az a)K_1^{n_{ijk}} \;.
\end{align}
\end{subequations}
Now, $X^i_{jk}\az a$ and $K_2^{-1}X^i_{jk}\az a$ are elements of
$\A(S^7_q)$ for all $a\in\A(S^4_q)$, hence are bounded since the
left regular representation of $\A(S^7_q)$ is bounded.
The restriction of $K_1$ to $\M_+$ is bounded too.
So, the left hand sides of (\ref{eq:Bcheck}), which are
the matrix elements of the operators in (\ref{eq:ltp}),
are bounded and this concludes the proof.
\end{prova}

For arbitrary (but fixed) $\lambda,\mu,\delta\in\C$ we define
\begin{equation}\label{eq:Dirac}
D_+:=\lambda X^1+\mu X^2+\delta X^3\;,\qquad
D_-:=(D_+)^*\;.
\end{equation}
The operator
$$
D=\ma{0 & D_- \\ D_+ & 0}\;,
$$
symmetric on $\M_+\oplus\M_-$ and odd with respect to the grading,
is our candidate for a Dirac operator. By Lemma \ref{lemma:Xbound},
the operators
$$
K_2^{-1}\bigl(D_+a -(K_2^{-1}\az a)D_+\bigr)\;,\qquad
\bigl(D_+a -(K_2^{-1}\az a)D_+\bigr)K_2^{-1}
$$
are bounded on $\M_+$. From this we deduce that the operator
\begin{align*}
& K_2^{-1}\bigl(D_-b-(K_2^{-1}\az b)D_-\bigr)=
K_2^{-1}\bigl(a^*D_- -D_-(K_2\az a^*)\bigr) \\
& \quad =K_2^{-1}\bigl(a^*D_- -D_-(K_2^{-1}\az a)^*\bigr)=
\bigl\{\bigl(D_+a -(K_2^{-1}\az a)D_+\bigr)K_2^{-1}\bigr\}^*
\end{align*}
is bounded on $\M_-$ for all $b\in\A$ (we called $a=-K_2\az b^*$).
Now $K_2$ does not map $\M_\pm$ into itself, so we cannot take
$K_2$ as twist operator. But if we define
\begin{equation}\label{eq:K}
K:=K_2\,\textrm{\footnotesize$\left(\!\!\begin{array}{cc}
q & 0 \\ 0 & 1\end{array}\!\!\right)$}
\end{equation}
then $K$ maps $\M_\pm$ into itself, $K^{-1}aK=K_2^{-1}\az a$ and
the operator
$$
\de a=K^{-1}\bigl(Da-(K^{-1}aK)D\bigr)=
\textrm{\footnotesize$\left(\!\!\begin{array}{cc}
0 & q^{-1}K_2^{-1}\bigl(D_+a -(K_2^{-1}\az a)D_+\bigr) \\
K_2^{-1}\bigl(D_-b-(K_2^{-1}\az b)D_-\bigr) & 0
\end{array}\!\!\right)$}
$$
is bounded for all $a\in\A$.

Although $D$ is only a symmetric operator, any odd symmetric
densely defined operator has a canonical selfadjoint extension
associated with the grading. Indeed, let $W$ be the restriction of
$\gamma$ to the range of $\,D+i\,$; since $W(D+i)=-(D-i)W$ the map
$W$ is the Cayley transform of $D$. Selfadjoint extensions of $D$
are in bijection with unitary extensions of $W$, but the grading
$\gamma$ \emph{is} a unitary extension of $W$, and this provides
the canonical extension that we need.

For any choice of the parameters $\lambda,\mu,\delta$ in Equation
(\ref{eq:Dirac}), the data $(\A(S^4_q),\HH,D,K)$ satisfies all the
axioms of an even twisted spectral triple, except for the compact
resolvent condition. Our hope is that it is possible to tune 
these three parameters so as to obtain an operator $D$ with compact
resolvent. To compute the spectrum of $D$ seems quite problematic,
and is postponed to future works.

\section{The real structure}\label{sec:sei}
The next step is to define the real structure.
Let $\kappa^{\frac{1}{2}}(a)=K_1^4K_2^3\az a\za K_1^4$ be the square
root of the modular automorphism (\ref{eq:mod}), and call
$T$ the antilinear operator on $\A(S^7_q)^2$ given by
\begin{equation}\label{eq:Treal}
T(v_1,v_2):=(\kappa^{\frac{1}{2}}(v_2^*),-\kappa^{\frac{1}{2}}(v_1^*)\bigr) \;.
\end{equation}
Since $\A(S^7_q)$ is an $U_q(su(2))$-module $*$-algebra,
we have $h\az v_i^*=\{S(h)^*\az v_i\}^*$. Moreover
$$
K_1\az\kappa^{\frac{1}{2}}(a)=\kappa^{\frac{1}{2}}(K_1\az a)\;,\qquad
E_1\az\kappa^{\frac{1}{2}}(a)=q^{-1}\kappa^{\frac{1}{2}}(E_1\az a)\;,\qquad
F_1\az\kappa^{\frac{1}{2}}(a)=q\kappa^{\frac{1}{2}}(F_1\az a)\;.
$$
Using these properties one checks that the operator $J_+$ (resp.~$J_-$),
defined by
$$
J_+(v_1,v_2)=T(qv_1,q^{-1}v_2)\;,\qquad
J_-(w_1,w_2)=T(q^{-1}w_1,qw_2)\;,
$$
maps $\M_+$ (resp.~$\M_-$) into itself.

Since $\kappa^{\frac{1}{2}}\circ *\circ\kappa^{\frac{1}{2}}\circ *=id$,
trivially $J_\pm^2=-1$. Furthermore, using the property (\ref{eq:twist}) and
the invariance of the Haar functional one proves that
\begin{align*}
\inner{Tv,Tw}&=
\sum\nolimits_i\varphi\bigl(\kappa^{\frac{1}{2}}(v_i^*)^*\,\kappa^{\frac{1}{2}}(w_i^*)\bigr)=
\sum\nolimits_i\varphi\bigl(\kappa^{-\frac{1}{2}}(v_i)\,\kappa^{\frac{1}{2}}(w_i^*)\bigr) \\ &=
\sum\nolimits_i\varphi\bigl(\kappa^{\frac{1}{2}}(w_i^*v_i)\bigr)=
\sum\nolimits_i\varphi(w_i^*v_i)=\inner{w,v}
\end{align*}
so that $T$ is an isometry and $J_\pm$ are bounded antilinear
operators and extend to $\HH_\pm$.

\begin{lemma}
If the parameters in Equation (\ref{eq:Dirac}) satisty
$\lambda=\bar{\mu}$ and $\delta=0$,
the operator $J=J_+\oplus J_-$ is a real structure for the
data $(\A(S^4_q),\HH,D,K)$.
\end{lemma}
\begin{prova}
The operator $J$ satisfies two of the three conditions in (\ref{eq:J}),
namely $J^2=-1$ and $J\gamma=\gamma J$. Since both $a$ and $\de a$
are operators of left multiplication for matrices with
entries in $\A(S^7_q)$ (cf.~the proof of Lemma \ref{lemma:Xbound}),
since
$$
JbJ^{-1}v=v\,\kappa^{\frac{1}{2}}(b^*)
$$
and since left and right multiplication commute, also the conditions
in Equation (\ref{eq:real}) are satisfied.
To prove that $J$ is a real structure we still have to prove that $DJ=JD$.
Since $J_-D_--D_-J_-=(D_+J_+-J_+D_+)^*$, it is sufficient to prove that
$D_+J_+=J_+D_+$. The first step is to observe that
$$
Th\textrm{\footnotesize$\left(\!\!\begin{array}{cc}
c_{11} & c_{12} \\ c_{21} & c_{22}\end{array}\!\!\right)$}T^{-1}=
K_1^4K_2^3S(h)^*K_1^{-4}K_2^{-3}
\textrm{\footnotesize$\left(\!\!\begin{array}{cc}
\bar{c}_{22} & -\bar{c}_{21} \\ -\bar{c}_{12} & \bar{c}_{11}
\end{array}\!\!\right)$}
$$
for all $h\in U_q(so(5))$ and $c_{ij}\in\C$.
With this formula we compute
$$
JX^1J^{-1}=X^2 \;,\qquad
JX^2J^{-1}=X^1 \;,\qquad
JX^3J^{-1}=K_2^2\textrm{\footnotesize$\left(\!\!\begin{array}{cc}
           q^{-1} & 0 \\ 0 & q\end{array}\!\!\right)$}X^3 \;.
$$
Thus if $\lambda=\bar{\mu}$ and $\delta=0$,
$D_+=\lambda X^1+\bar{\lambda}X^2$ and $JD_+=D_+J$.
This concludes the proof.
\end{prova}



\end{document}